\documentclass[a4paper,11pt,oneside]{amsart}
\usepackage{amssymb, color, graphics, colordvi,euscript}
\usepackage{amscd}
\usepackage{graphicx}
\usepackage{epsfig}
\usepackage{xypic}

\newcommand{\bb}{\mathbb}
\newcommand{\cc}{\bb C}
\newcommand{\rr}{\bb R}
\newcommand{\R}{\rr}

\newcommand{\pp}{\bb P}
\renewcommand{\P}{\mathbb{P}}

\newcommand{\Diff}{\operatorname{Diff}}
\newcommand{\Aut}{\operatorname{Aut}}

\newcommand{\SX}{\mathcal{X}}

\let\ra\rightarrow

\let\iso\cong

\newcommand{\so}{\operatorname{SO}}

\renewcommand{\phi}{\varphi}
\renewcommand{\tilde}{\widetilde}

\numberwithin{equation}{section}%
\newtheorem{theo}[equation]{Theorem}%

\newtheorem{lem}[equation]{Lemma}

\newtheorem{cor}[equation]{Corollary}

\newtheorem{con}[equation]{Conjecture}

\theoremstyle{remark}

\newtheorem{que}[equation]{Question}
\newtheorem{ack}{Acknowledgements}

\begin{document}

\title[Automorphisms of real rational surfaces]{The group of
automorphisms of a real rational surface is $n$-transitive}

\dedicatory{To Joost van Hamel in memoriam}

\author{ Johannes Huisman \and Fr\'ed\'eric Mangolte }
\address{Johannes Huisman, D\'epartement de Math\'ematiques,
 Laboratoire CNRS UMR 6205, Universit\'e de Bretagne Occidentale, 6,
 avenue Victor Le Gorgeu, CS~93837, 29238 Brest cedex 3, France.
 Tel.~+33 2 98 01 61 98, Fax~+33 2 98 01 67 90}
\email{johannes.huisman@univ-brest.fr}
\urladdr{http://pageperso.univ-brest.fr/$\sim$huisman}
\address{Fr\'ed\'eric Mangolte, Laboratoire de Math\'ematiques,
 Universit\'e de Savoie, 73376 Le Bourget du Lac Cedex, France,
 Phone: +33 (0)4 79 75 86 60, Fax: +33 (0)4 79 75 81 42}
\email{mangolte@univ-savoie.fr}
\urladdr{http://www.lama.univ-savoie.fr/$\sim$mangolte}

\thanks{The research of the second author was partially supported by the  
ANR grant "JCLAMA" of the french "Agence Nationale de la Recherche".}

\date{}

\begin{abstract}
Let~$X$ be a rational nonsingular compact connected real algebraic
surface.  Denote by~$\Aut(X)$ the group of real algebraic automorphisms of~$X$.  We
show that the group~$\Aut(X)$ acts $n$-transitively on~$X$, for all
natural integers~$n$.

 As an application we give a new and simpler proof of the fact that
 two rational nonsingular compact connected real algebraic surfaces
 are isomorphic if and only if they are homeomorphic
 as topological surfaces.
\end{abstract}

\maketitle


\begin{quote}\small
\end{quote}

\section{Introduction}\label{sec:intro}

Let~$X$ be a nonsingular compact connected real algebraic manifold,
i.e., $X$ is a compact connected submanifold of~$\R^n$ defined by real
polynomial equations, where~$n$ is some natural integer. We study the
group of algebraic automorphisms of~$X$. Let us make precise what we
mean by an algebraic automorphism.

Let $X$~and $Y$ be real algebraic submanifolds of $\R^n$ and $\R^m$,
respectively.  An \emph{algebraic map}~$\phi$ of~$X$ into~$Y$ is a
map of the form
\begin{equation}\label{eqalgmap}
\phi(x)=\left(\frac{p_1(x)}{q_1(x)},\ldots,\frac{p_m(x)}{q_m(x)}\right)
\end{equation}
where~$p_1,\ldots,p_m,q_1,\ldots,q_m$ are real polynomials in the
variables~$x_1,\ldots,x_n$, with~$q_i(x) \neq0$ for any~$x \in X$ and
any~$i$.  An algebraic map from~$X$ into~$Y$ is also called a
\emph{regular map}~\cite{BCR}.  Note that an algebraic map is
necessarily of class~$C^\infty$.  An algebraic map $\phi \colon X \to
Y$ is an \emph{algebraic isomorphism}, or \emph{isomorphism} for
short, if $\phi$ is algebraic, bijective and if~$\phi^{-1}$ is
algebraic.  An algebraic isomorphism from $X$~into $Y$ is also called
a \emph{biregular map}~\cite{BCR}.  Note that an algebraic isomorphism
is a diffeomorphism of class~$C^\infty$. An algebraic isomorphism
from~$X$ into itself is called an \emph{algebraic automorphism}
of~$X$, or \emph{automorphism} of~$X$ for short.  We denote
by~$\Aut(X)$ the group of automorphism of~$X$.

For a general real algebraic manifold, the group~$\Aut(X)$ tends to be
rather small. For example, if~$X$ admits a complexification~$\SX$ that
is of general type then~$\Aut(X)$ is finite.  Indeed, any automorphism
of~$X$ is the restriction to~$X$ of a birational automorphism
of~$\SX$. The group of birational automorphisms of~$\SX$ is known to
be finite~\cite{Matsumura}.  Therefore, $\Aut(X)$ is finite for such
real algebraic manifolds.

In the current paper, we study the group~$\Aut(X)$ when~$X$ is a
compact connected real algebraic surface, i.e., a compact connected
real algebraic manifold of dimension~$2$. By what has been said above, the
group of automorphisms of such a surface is most
interesting when the Kodaira dimension of~$X$ is equal to $-\infty$,
and, in particular, when~$X$ is geometrically rational.  By a
result of Comessatti, a connected geometrically rational real surface
is rational (see Theorem~IV of~\cite{Co12} and the remarks
thereafter, or \cite[Corollary~VI.6.5]{Si}). Therefore, we will
concentrate our attention to the group~$\Aut(X)$ when $X$ is a
rational compact connected real algebraic surface.

Recall that a real algebraic surface~$X$ is \emph{rational} if there
are a nonempty Zariski open subset $U$ of~$\R^2$, and a nonempty
Zariski open subset~$V$ of~$X$, such that $U$~and $V$ are isomorphic
real algebraic varieties, in the sens above.  In particular, this
means that $X$ contains a nonempty Zariski open subset $V$ that admits
a parametrization by real rational functions in two variables.

Examples of rational real algebraic surfaces are the following:
\begin{itemize}
\item the unit sphere~$S^2$ defined by the equation $x^2+y^2+z^2=1$
in~$\R^3$,
\item the real algebraic torus~$S^1\times S^1$, where $S^1$ is the
 unit circle defined by the equation~$x^2+y^2=1$ in~$\R^2$, and
\item any real algebraic surface obtained from one of the above ones
 by repeatedly blowing up a point.
\end{itemize}
This is a complete list of rational real algebraic surfaces, as was
probably known already to Comessatti. A modern proof may use the
Minimal Model Program for real algebraic surfaces~\cite{Ko1,Ko2}
(cf.~\cite[Theorem~3.1]{BH07}).  For example, the real projective
plane~$\P^2(\R)$|of which an explicit realization as a rational real
algebraic surface can be found in~\cite[Theorem~3.4.4]{BCR}|is
isomorphic to the real algebraic surface obtained from $S^2$ by
blowing up $1$ point.

The following conjecture has attracted our attention.

\begin{con}[\protect{\cite[Conjecture~1.4]{BH07}}]\label{cocon}
 Let~$X$ be a rational nonsingular compact connected real algebraic
 surface. Let~$n$ be a natural integer.  Then the
 group~$\Aut(X)$ acts $n$-transitively on~$X$.
\end{con}

The conjecture seems known to be true only in the case
when~$X$ is isomorphic to~$S^1\times S^1$:

\begin{theo}[\protect{\cite[Theorem~1.3]{BH07}}]
\label{ths1xs1}
The group~$\Aut(S^1\times S^1)$ acts $n$-transi\-tively
on~$S^1\times S^1$, for any natural integer~$n$.\qed
\end{theo}

The object of the paper is to prove Conjecture~\ref{cocon}:

\begin{theo}\label{theo:main}
 The group $\Aut(X)$ acts $n$-transitively on $X$, whenever $X$
 is a rational nonsingular compact connected real algebraic surface,
 and $n$ is a natural integer.
\end{theo}

Our proof goes as follows. We first prove $n$-transitivity
of~$\Aut(S^2)$ (see Theorem~\ref{theo:s2}). For this, we need a
large class of automorphisms of~$S^2$.
Lemma~\ref{lem:phi} constructs such a large class. Once
$n$-transitivity of~$\Aut(S^2)$ is established, we prove
$n$-transitivity of~$\Aut(X)$, for any other rational
surface~$X$, by the following argument.

If~$X$ is isomorphic to~$S^1\times S^1$ then the $n$-transitivity has
been proved in~\cite[Theorem~1.3]{BH07}.  Therefore, we may assume
that~$X$ is not isomorphic to~$S^1\times S^1$. We prove that~$X$ is
isomorphic to a blowing-up of~$S^2$ in $m$ distinct points, for some
natural integer~$m$ (see Theorem~\ref{thmmps} for a precise
statement). The $n$-transitivity of~$\Aut(X)$ will then follow from
the $(m+n)$-transitivity of~$\Aut(S^2)$.

Theorem~\ref{theo:main} shows that the group of automorphisms of a
rational real algebraic surface is big.  It would, therefore, be
particularly interesting to study the dynamics of automorphisms of
rational real surfaces, as is done for K3-surfaces in~\cite{Ca01}, for
example.

Using the results of the current paper, we were able, in a forthcoming
paper~\cite{HM08}, to generalize Theorem~\ref{theo:main} and prove
$n$-transitivity of~$\Aut(X)$ for curvilinear infinitely near points on a rational
surface~$X$.

We also pass to the reader the following interesting question of the
referee.

\begin{que}
Let $X$ be a rational nonsingular compact connected real algebraic
surface.  Is the subgroup~$\Aut(X)$ dense in the group~$\Diff(X)$ of
all $C^\infty$ diffeomorphisms of~$X$ into itself?
\end{que}

This question is studied in the forthcoming paper \cite{KM08}. 

As an application of Theorem~\ref{theo:main}, we present in
Section~\ref{seBH} a simplified proof of the following result.

\begin{theo}[\protect{\cite[Theorem~1.2]{BH07}}]\label{thratmod}
 Let~$X$ and $Y$ be rational nonsingular compact connected real
 algebraic surfaces.  Then the following statements are equivalent.
\begin{enumerate}
\item The real algebraic surfaces $X$ and $Y$ are isomorphic.
\item The topological surfaces $X$ and $Y$ are homeomorphic.
\end{enumerate}
\end{theo}

\begin{ack}
We want to thank J.~Koll\'ar for useful discussions, in
particular about terminology.  We are also grateful to the referee
for helpful remarks that allowed us to improve the exposition.
\end{ack}

\section{$n$-Transitivity of $\protect{\Aut}(S^2)$}

We need to slightly extend the notion of an algebraic map between real
algebraic manifolds.  Let $X$~and $Y$ be real algebraic submanifolds
of~$\R^n$ and $\R^m$, respectively.  Let~$A$ be any subset of~$X$. An
\emph{algebraic map} from $A$~into $Y$ is a map~$\phi$ as
in~(\ref{eqalgmap}), where~$p_1,\ldots,p_m,q_1,\ldots,q_m$ are real
polynomials in the variables~$x_1,\ldots,x_n$, with~$q_i(x) \neq0$ for
any~$x \in A$ and any~$i$.  To put it otherwise, a map~$\phi$ from~$A$
into $Y$ is \emph{algebraic} if there is a Zariski open subset~$U$
of~$X$ containing~$A$ such that~$\phi$ is the restriction of an
algebraic map from~$U$ into~$Y$.

We will consider algebraic maps from a subset~$A$ of~$X$ into~$Y$, in
the special case where $X$ is isomorphic to the
real algebraic line~$\R$, the subset~$A$ of~$X$ is a closed
interval, and $Y$ is isomorphic to the real
algebraic group~$\so_2(\R)$.

Denote by $S^2$ the 2-dimensional sphere defined in $\rr^3$ by the equation
$$
x^2+y^2+z^2=1.
$$

\begin{lem}\label{lem:phi}
 Let $L$ be a line through the origin of $\rr^3$ and denote by $I
 \subset L$ the closed interval whose boundary is $L \cap S^2$.
 Denote by $L^\perp$ the plane orthogonal to $L$ containing the
 origin. Let $f \colon I \to \so(L^\perp)$ be an
 algebraic map. Define $\phi_f \colon S^2 \to S^2$ by
$$
\phi_f (z,x) = (f(x)z,x)
$$
where $(z,x) \in (L^\perp \oplus L) \cap S^2$.
Then $\phi_f$ is an automorphism of~$S^2$.
\end{lem}

\begin{proof}
 Identifying $\R^2$ with~$\cc$, we may assume that $S^2 \subset \cc
 \times \rr$ is given by the equation~$\vert z \vert ^2 + x^2 =1$,
 and that $L$ is the line $\{0\}\times \rr$. Then $L^\perp = \cc
 \times \{0\}$ and $\so(L^\perp) = S^1$.  It is clear that the map
 $\phi_f$ is an algebraic map from $S^2$ into itself. Let
 $\overline{f}$ be the complex conjugate of $f$, i.e. $\forall x \in
 I,\ \overline{f}(x) = \overline{f(x)}$. We have $\phi_{\overline{f}}
 \circ \phi_f = \phi_f \circ \phi_{\overline{f}} = id$. Therefore
 $\phi_f$ is an automorphism of~$S^2$.
\end{proof}

\begin{lem}\label{lem:interpolation}
 Let $x_1,\dots,x_n$ be $n$ distinct points of the closed
 interval~$[-1,1]$, and let $\alpha_1,\dots,\alpha_n$ be elements of
 $\so_2(\rr)$.  Then there is an algebraic map $f
 \colon [-1,1] \to \so_2(\rr)$ such that $f(x_j) =
 \alpha_j$ for $j=1,\dots ,n$.
\end{lem}

\begin{proof}
 Since $\so_2(\rr)$ is isomorphic to the unit circle
 $S^1$, it suffices to prove the statement for $S^1$ instead of
 $\so_2(\rr)$. Let~$P$ be a point of~$S^1$ distinct
 from~$\alpha_1,\ldots,\alpha_n$.  Since $S^1\setminus \{P\}$ is
 isomorphic to $\rr$, it suffices, finally, to prove
 the statement for $\rr$ instead of $\so_2(\rr)$. The
 latter statement is an easy consequence of Lagrange polynomial
 interpolation.
\end{proof}

\begin{theo}\label{theo:s2}
 Let $n$ be a natural integer. The group $\Aut(S^2)$ acts
 $n$-transitively on $S^2$.
\end{theo}

\begin{proof}
We will need the following terminology.  Let $W$ be a point of
$S^2$, let $L$ be the line in $\rr^3$ passing through $W$ and the
origin. The intersection of $S^2$ with any plane in~$\R^3$ that is
orthogonal to $L$ is called a {\em parallel of~$S^2$ with respect to
  $W$}.

Let $P_1,\dots,P_n$ be $n$ distinct points of~$S^2$, and let
$Q_1,\dots,Q_n$ be $n$ distinct points of~$S^2$.  We need to show that
there is an automorphism~$\phi$ of~$S^2$
such that~$\phi(P_j)=Q_j$, for all~$j$.

 Up to a projective linear automorphism of $\pp^3(\rr)$ fixing $S^2$,
 we may assume that all the points $P_1,\dots,P_n$ and
 $Q_1,\dots,Q_n$ are in a sufficiently small neighborhood of the
 north pole~$N = (0,0,1)$ of~$S^2$. Indeed, we may first assume that none
 of these points is contained in a small spherical disk $D$ centered at
 $N$. Then the images of the points by the inversion with respect to
 the boundary of $D$ are all contained in $D$.

 We can choose two points $W$ and $W'$ of $S^2$ in the $xy$-plane
 such that the angle $WOW'$ is equal to $\pi/2$ and such that the
 following property holds. Any parallel with respect to $W$ contains
 at most one of the points $P_1,\dots,P_n$, and any parallel with
 respect to $W'$ contains at most one of $Q_1,\dots,Q_n$. Denote by
 $\Gamma_j$ the parallel with respect to $W$ that contains $P_j$, and
 by $\Gamma_j'$ the one with respect to $W'$ that contains $Q_j$.

 Since the disk $D$ has been chosen sufficiently small, $\Gamma_j
 \cap \Gamma_j'$ is nonempty for all $j=1,\dots, n$. Let $R_j$ be one
 of the intersection points of $\Gamma_j$ and $\Gamma_j'$ (see
 Figure~\ref{fisphere}).  It is now sufficient to show that there is
 an automorphism $\varphi$ of $S^2$ such that
 $\varphi(P_j)= R_j$.

\begin{figure}
\centering\leavevmode
 \epsfbox{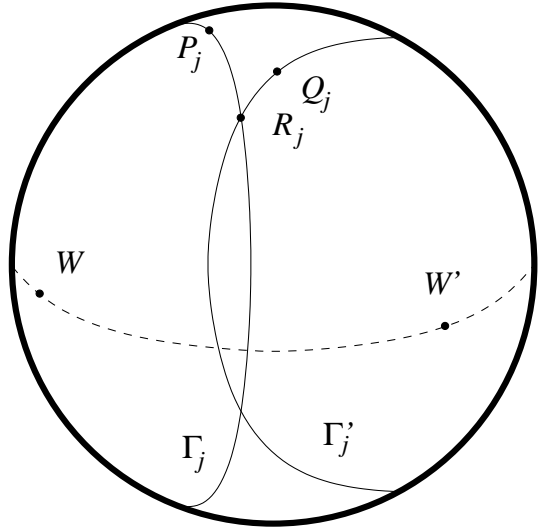} \caption{The sphere $S^2$ with the parallels $ 
\Gamma_j$ and $\Gamma_j'$.}
 \label{fisphere}
\end{figure}

 Let again $L$ be the line in $\rr^3$ passing through $W$ and the
 origin. Denote by $I \subset L$ the closed interval whose boundary
 is $L \cap S^2$. Let $x_j$ be the unique element of $I$ such that
 $\Gamma_j = (x_j + L^\perp) \cap S^2$. Let $\alpha_j \in
 \so(L^\perp)$ be such that $\alpha_j (P_j - x_j) = R_j
 - x_j$. According to Lemma~\ref{lem:interpolation}, there is an
 algebraic map $f \colon I \to \so(L^\perp)$ such that
 $f(x_j) = \alpha_j$. Let $\phi := \phi_f$ as in Lemma~\ref{lem:phi}.
 By construction, $\phi(P_j) = R_j$, for all $j=1,\dots,n$.
\end{proof}

\section{$n$-Transitivity of~$\Aut(X)$}

\begin{theo}\label{thmmps}
 Let~$X$ be a rational nonsingular compact connected real algebraic
 surface and let~$S$ be a finite subset of~$X$. Then,
\begin{enumerate}
\item $X$ is either isomorphic to~$S^1\times S^1$, or
\item there are distinct points~$R_1,\ldots,R_m$ of~$S^2$ and a finite
subset~$S'$ of~$S^2$ such that
\begin{enumerate}
\item $R_1,\ldots,R_m\not\in S'$, and
\item there is an isomorphism~$\phi\colon X\ra
B_{R_1,\ldots,R_m}(S^2)$ such that~$\phi(S)=S'$.
\end{enumerate}
\end{enumerate}
\end{theo}

\begin{proof}
By what has been said in the introduction, $X$ is either isomorphic
to~$S^1\times S^1$, in which case there is nothing to prove, or~$X$
is isomorphic to a real algebraic surface obtained from~$S^2$ by
successive blow-up. Therefore, we may assume that there is a sequence
$$
\xymatrix{X=X_m\ar[r]^{f_m}&X_{m-1}\ar[r]^{f_{m-1}}&
 \cdots\ar[r]^{f_1}&X_0=S^2},
$$ 
where $f_i$ is the blow-up of~$X_{i-1}$ at a point~$R_i$ of~$X_{i-1}$.

Let~$\tilde{S}$ be the union of~$S$ and the set of
centers~$R_1,\ldots,R_m$. Since the elements of~$\tilde{S}$ can be
seen as infinitely near points of~$S^2$, there is a natural partial
ordering on~$\tilde{S}$. The partially ordered set~$\tilde{S}$ is a
finite forest with respect to that ordering.  

The statement that we need to prove is that there is a sequence of
blow-ups as above such that all trees of the corresponding forest
have height~$0$. We prove that statement by induction on the sum~$h$
of heights of the trees of the forest~$\tilde{S}$. If~$h=0$ there is nothing
to prove.  Suppose, therefore, that~$h\neq0$.  We may then assume,
renumbering the $R_i$ if necessary, that either~$R_2\leq R_1$ or
that a point~$P\in S$ is mapped onto~$R_1$ by the
composition~$f_n\circ\cdots\circ f_1$.

As we have mentioned in the introduction, the real algebraic surface
obtained from~$S^2$ by blowing up at~$R_1$ is isomorphic to the real
projective plane~$\P^2(\R)$. Moreover, the exceptional divisor
in~$\P^2(\R)$ is a real projective line~$L$. We
identify~$B_{R_1}(S^2)$ with~$\P^2(\R)$.  Choose a real projective
line~$L'$ in~$\P^2(\R)$ such that no element
of~$\tilde{S}\setminus\{R_1\}$ is mapped into~$L'$ by a suitable
composition of some of the maps $f_2,\ldots, f_m$. Since the group of
linear automorphisms of~$\P^2(\R)$ acts transitively on the set of
projective lines, the line~$L'$ is an exceptional divisor for a
blow-up~$f_1'\colon \P^2(\R)\ra S^2$ at a point~$R_1'$ of~$S^2$.  It
is clear that the sum of heights of the trees of the corresponding
forest is equal to~$h-1$. The statement of the theorem
follows by induction.
\end{proof}

\begin{cor}\label{commp}
 Let~$X$ be a rational nonsingular compact connected real algebraic
 surface. Then,
\begin{enumerate}
\item $X$ is either isomorphic to~$S^1\times S^1$, or
\item there are distinct points~$R_1,\ldots,R_m$ of~$S^2$ such
that~$X$ is isomorphic to the real algebraic surface obtained
from~$S^2$ by blowing up the points~$R_1,\ldots,R_m$.\qed
\end{enumerate}
\end{cor}

\begin{proof}[Proof of Theorem~\ref{theo:main}]
Let $X$ be a rational surface and let~$(P_1,\ldots,P_n)$ and
$(Q_1,\ldots,Q_n)$ by two $n$-tuples of disctinct points of~$X$. By
Theorem~\ref{thmmps}, $X$ is either isomorphic to~$S^1\times S^1$ or
to the blow-up of~$S^2$ at a finite number of distinct
points~$R_1,\ldots,R_m$.  If $X$ is isomorphic to~$S^1\times S^1$
then~$\Aut(X)$ acts $n$-transitively by~\cite[Theorem~1.3]{BH07}.
Therefore, we may assume that~$X$ is the
blow-up~$B_{R_1,\dots,R_m}(S^2)$ of~$S^2$ at~$R_1,\ldots,R_m$.
Moreover, we may assume that the points
$P_1,\ldots,P_n,Q_1,\ldots,Q_n$ do not belong to any of the
exceptional divisors. This means that these points are elements
of~$S^2$, and that, $(P_1,\ldots,P_n)$ and
$(Q_1,\ldots,Q_n)$ are two $n$-tuples of distinct points of~$S^2$. It
follows that $(R_1,\ldots,R_m,P_1,\ldots,P_n)$ and
$(R_1,\ldots,R_m,Q_1,\ldots,Q_n)$ are two $(m+n)$-tuples of distinct
points of~$S^2$.  By Theorem~\ref{theo:s2}, there is an automorphism
$\psi$ of~$S^2$ such that~$\psi(R_i)=R_i$, for all~$i$,
and~$\psi(P_j)=Q_j$, for all~$j$.  The induced automorphism~$\phi$
of~$X$ has the property that~$\psi(P_j)=Q_j$, for all~$j$.
\end{proof}

\section{Classification of rational real algebraic
surfaces}
\label{seBH}

\begin{proof}[Proof of Theorem~\ref{thratmod}]
Let~$X$ and $Y$ be a rational nonsingular compact connected real
algebraic surfaces.  Of course, if $X$~and $Y$ are isomorphic then
$X$~and $Y$ are homeomorphic. In order to prove the converse,
suppose that $X$~and $Y$ are homeomorphic. We show that there is an
isomorphism from~$X$ onto~$Y$.

By Corollary~\ref{commp}, we may assume that $X$~and $Y$ are not
homeomorphic to~$S^1\times S^1$. Then, again by
Corollary~\ref{commp}, $X$~and $Y$ are both isomorphic to a real
algebraic surface obtained from~$S^2$ by blowing up a finite number
of distinct points. Hence, there are distinct points $P_1,\ldots,P_n$
of~$S^2$ and distinct points $Q_1,\ldots,Q_m$ of~$S^2$ such that
$$
X\iso B_{P_1,\ldots,P_n}(S^2) \quad\text{and}\quad Y\iso
B_{Q_1,\ldots,Q_m}(S^2).
$$
Since~$X$ and $Y$ are homeomorphic, $m=n$. By Theorem~\ref{theo:s2},
there is an automorphism~$\phi$ from~$S^2$ into $S^2$ such that
$\phi(P_i)=Q_i$ for all~$i$. It follows that~$\phi$ induces an
algebraic isomorphism from~$X$ onto~$Y$.
\end{proof}


\end{document}